\documentclass[12, reqno]{amsart}
\usepackage{amsmath}
\usepackage{amscd,amsthm,amsfonts,amsopn,amssymb,verbatim}
\usepackage{mathrsfs}
\parskip =\medskipamount
\numberwithin{equation}{section}
\newtheorem{theorem}{Theorem} 

\newtheorem{proposition}[theorem]{Proposition}

\theoremstyle{remark}

\DeclareMathOperator{\supp}{supp\,}

\def\be{\begin{equation}}
\def\ee{\end{equation}}
\def\vp{\varphi}
\def\ve{\varepsilon}

\def\nint{\mathop{\diagup\kern-13.0pt\int}}

\allowdisplaybreaks 
\begin{document} 
\large
\setlength{\parskip}{2pt}

\title{On random walks in large compact Lie groups}
\date{\today}
\author{Jean Bourgain}
\address{(J. Bourgain) Institute for Advanced Study, Princeton, NJ 08540}
\email{bourgain@math.ias.edu}
\thanks{This work was partially supported by NSF grants DMS-1301619}
\maketitle 

\section
{Introduction}

In order to put the problem considered in this Note in perspective, we first recall some other relatively recent results
around spectral gaps and generation in Lie groups.

It was shown in \cite {B-G1} (resp. \cite {B-G2}) that if $\Lambda$ is a symmetric finite subset of $SU(2)$ \big(resp.
$SU(d)$\big) consisting of algebraic elements, such that the countable group $\Gamma =\langle\Lambda\rangle$
generated by $\Lambda$ is dense, then the corresponding averaging operators
$$
Tf=\frac 1{|\Lambda|}\sum_{g\in\Lambda} f\circ g\eqno{(1.1)}
$$
acting on $L^2(G)$, has a uniform spectral gap (only depending on $\Lambda$).
This result was generalized in $[dS-B]$ to simple compact Lie groups.

It is not known if the assumption for $\Lambda$ to be algebraic is needed, and one may conjecture that it is not.
Short of providing uniform spectral gaps, Varju \cite {V} established the following property which is the most relevant
statement for what follows.

\begin{proposition}
Let $G$ be a compact Lie group with semisimple connected component.
Let $\mu$ be a probability measure on $G$ such that {\rm supp}\,$(\tilde\mu *\mu)$, $\tilde\mu$ defined by $\int f(x)d\tilde
\mu(x) =\int f(x^{-1})d\mu(x)$, generates a dense subgroup of $G$.
Then there is a constant $c>0$ depending only on $\mu$ such that the following holds.

Let $\vp \in {\rm Lip } (G), \Vert\vp\Vert_2 =1$ and $\int_G\vp=0$.
Then
$$
\Big\Vert\int \vp (h^{-1} g) d\mu(h)\Big\Vert_2 < 1-c \, \log ^{-A} (1+\Vert\vp \Vert_{\rm Lip})\eqno{(1.2)}
$$
with $A$ depending on $G$.
\end{proposition}

Using (1.2) and decomposition of the regular representation of $G$ in irreducibles (though this may be avoided), one
deduces easily from (1.2) that it takes time at most $O(\log ^A \frac 1\ve)$ as $\ve\to 0$ for the random walk governed by
$\mu$ to produce an $\ve$-approximation of uniform measure on $G$.
Note that for $G=SU(d)$, this statement corresponds to the Solovay-Kitaev estimates on generation, cf. \cite{D-N}, which
in fact turns out to be equivalent.

Let us focus on $G=SO(d)$ or $SU(d)$.
While the exponent $A$ in (1.2) is a constant, the prefactor $c$ depends on $\mu$, hence on $G$, and seems to have
received little attention.
Basically our aim is to prove a lower bound on $c$ which is powerlike in $\frac 1d$ and without the need for uniform
spectral gaps (which may not be always available).
We focus on the following model problem brought
to the author's attention by T.~Spencer (who was motivated by issues in random matrix
theory that will not be pursued here).
The general setting is as follows (we consider the $SU(d)$-version).
Fix some probability measure $\eta$ on $SU(2)$ such that its support generates a dense group, i.e.
$\overline{\langle \supp \eta\rangle} = SU(2)$.
This measure $\eta$ may be Haar but could be taken discrete as well.
Identify $\{0, 1, \ldots, d-1\}$ with the cyclic group $\mathbb Z/d\mathbb Z$ and denote $\nu_{ij}$ the measure $\eta$
on $SU(2)$ acting on
the space $[e_i, e_j]$.
Consider the random walk on $SU(d)$ given by
$$
Tf(x) =\frac 1d \sum^{d-1}_{i=0} \int f(gx) \nu_{i, i+1}(dg).\eqno{(1.3)}
$$

How long does it take for this random walk to become an $\ve$-approximation of uniform measure on $G$, with special
emphasis on large $d$?
Thus this is a particular instance of the more general issue formulated in the title.
While we are unable to address the broader problem, specific cases such as (1.3) 
may be analyzed in a satisfactory way (based partly on
arguments that are also relevant to the general setting).

We prove

\begin{proposition}
In the above setting, $\ve$-approximation of the uniform measure is achieved in time $C(d\log\frac 1\ve)^C$, with $C$ a
constant independent of $d$.
\end{proposition}

\noindent
{\bf Comment}

If $\eta$ is taken to be a uniform measure on $SU(2)$, better results are available, exploiting Hurwitz' construction of
Haar measure (see \cite{D-SC}, section 2).
In this situation, the operator $T$ displays in fact a uniform spectral gap and the power of $\log\frac 1\ve$ can be
taken to be one (cf. \cite{D-SC}, Theorem 1).
Our interest in this presentation is a more robust approach however.

Basically, one could expect a more general phenomenon (though some additional assumptions are clearly needed).
In some sense, it would give a continuous version of the conjecture of Babai and Seress \cite{B-S} predicting
poly-logarithmic diameter for the family of non-Abelian finite simple groups (independently of the choice of generators).
Important progress in this direction for the symmetric group appears in \cite {H-S}.

Independently of Spencer's question, related spectral gap and mixing time issues for specific random walks in large (not
necessarily compact) linear groups appear in the theory of Anderson localization for `quasi-one-dimensional' methods in
Math Phys.

Consider the strip $\mathbb Z\times\mathbb Z/d\mathbb Z$ and a random Schr\"odinger operator $\Delta +\lambda V$ with
$\Delta$ the usual lattice Laplacian on $\mathbb Z\times\mathbb Z/ d\mathbb Z$, $V$ a random potential and $\lambda>0$
the disorder.
This model is wellknown to exhibit pure point spectrum with so-called Anderson localization for the eigenfunctions.
The issue here is how the localization length (or equivalently, the Lyapounov exponents in the transfer matrix approach)
depend on $d$ when $d\to\infty$.

The classical approach based on Furstenberg's random matrix product theory (acting on extension powers of $\mathbb R^d$),
cf. \cite {B-L}, is not quantitative and sheds no light on the role of $d$.
In fact, the first explicit lower bound on Lyapounov exponents seems to appear in \cite {B} (using different techniques
based on Green's function analysis), with, roughly speaking exponential dependence on $d$ (while the `true' behaviour is
believed to be rather of the form $d^{-C}$).
Clearly understanding the mixing time for the random walk in the symplectic group $Sp(2d)$ associated to the 
transfer matrix is crucial.
Note that this group is non-compact, which is an added difficulty 
(for very small $\lambda$, depending on $d$, \cite {B-S} provides the precise asymptotic of the exponents, based on a
multi-dimensional extension of the Figotin-Pastur approach).

\section
{Some preliminary comments}

The proof of Proposition 1 in \cite{V} exploits the close relation between `generation' and `restricted spectral gaps'.
This point of view is also the key idea here in establishing

\medskip

\noindent
{\bf Proposition 1$'$}
{\sl Let $T$ be defined by (1.3).
Then there is the following estimate
$$
\Vert T f\Vert_2 <1- (cd)^{-C} \big(\log (1+\Vert f\Vert _{\rm Lip})^{-A}\big) \eqno{(2.1)}
$$
for $f\in {\rm Lip} (G)$. $\Vert f\Vert_2 =1, \int_G f=0$.}

\medskip

Here $C$ and $A$ are constants (denoted differently, because of their different appearance in the argument).
\medskip

Unlike in \cite{V}, we tried to avoid the use of representation theory.
The reason for this is the following.
If one relies on decomposition of the regular representation of $G$ in irreducibles and the Peter-Weyl theorem, one is
faced in the absence of a uniform spectral gap with convergence issues of the generalized Fourier expansion of functions on $G$ of given regularity.
Conversely, we also need to understand the regularity of matrix coefficients of the representations of increasing
dimension.
While these are classical issues, understanding the role of the dimension $d$ does not seem to have been addressed
explicitly.

\medskip

\section
{Proof of proposition 1$'$}

For simplicity, we take $\eta$ to be a uniform measure on $SU(2)$ and indicate the required modifications for the general case
in \S5.

According to (1.3), denote
$$
\nu =\frac 1d \sum_{i=0}^{d-1} \nu_{i, i+1}\eqno{(3.1)}
$$
Thus $\nu =\tilde\nu$ and $T$ is the corresponding averaging operator.

Let $f\in {\rm Lip} (G), \Vert f\Vert_2=1$ and $\int_G f=0$.
Assuming 
$$
\Big\Vert \int \tau_g f \nu(dg)\Big \Vert^2_2 = \Vert Tf\Vert^2_2 > 1-\ve\eqno{(3.2)}
$$
(denoting $\tau_g f(x) = f(gx)\big)$
our aim is to obtain a lower bound on $\ve$.

Clearly (3.2) implies that
$$
\Big\langle f, \int \tau_gf(\nu*\nu)(dg)\Big\rangle > 1-\ve
$$
and
$$
\int\Vert f-\tau_g f\Vert^2_2 (\nu*\nu) (dg)< 2\ve.\eqno{(3.3)}
$$

Fix $\ve_1>0$ to be specified later and denote $B_{\ve_1}$ an
$\ve_1$-neighborhood (for the operator norm) of $Id$ in $SU(d)$.
It is clear from (3.1) that $\nu(B_{\ve_1})\gtrsim \ve_1^4$ and hence (3.3) implies
$$
\int \Vert f-\tau_{g'g}f\Vert^2_2 \, \nu(dg)\lesssim \ve_1^{-4}\ve\eqno{(3.4)}
$$
for some $g'\in B_{\ve_1}$.
Next, partitioning $SU(2)$ in $\ve_1$-cells $\Omega_\alpha$ and denoting
$$
\Omega_{\alpha, i}=\{ g\in SU(d); g(e_j)=e_j \ \text { for } \ j\not\in\{i, i+1\}\ \text { and } \ g|_{[e_i, e_{i+1}]}
\in\Omega_\alpha\}
$$
observe that $\nu(\Omega _{\alpha, i})\geq \frac 1d \ve_1^4$ so that by (3.4)
$$
\nint_{\Omega_{\alpha, i}} \Vert f- \tau_{g'g} f\Vert^2_2 \, \nu (dg) \lesssim d\ve_1^{-8} \ve\ll 1.\eqno{(3.5)}
$$
Exploiting (3.5), it is clear that we may introduce a collection
$\mathcal G\subset SU(d)$ with the following properties
$$
\Vert f-\tau_g f\Vert_2 \lesssim \sqrt d \, \ve_1^{-4} \sqrt\ve \ \text { for } \ g\in\mathcal G. \eqno{(3.6)}
$$
and 

Given an element $\gamma\in SU(2)$ and $1\leq i<j\leq d$, denote
$\gamma_{ij}$ in $SU(d)$ the element defined by
$$
\begin{cases}
\gamma_{ij}(e_k)&= e_k \ \text { for } \ k\not \in \{i, j\}\\[5pt]
\gamma_{ij} \big|_{[e_i, e_j]}&=\gamma.
\end{cases}
\eqno{(3.7)}
$$
Then, for each $\gamma\in SU(2)$ and $1\leq i\leq d$, there is $g\in \mathcal G$ s.t.
$$
\Vert g -\gamma_{i, i+1}\Vert_2 < \ve_1.\eqno{(3.8)}
$$

At this point, we will invoke generation.
Since $\int_G f =0$,
$$
\int_{SU(d)} \Vert f-\tau_g f\Vert^2_2 dg =2
$$
and we take some $h_0\in SU(d)$ s.t.
$$
\Vert f-\tau_{h_0}f\Vert_2\geq \sqrt 2.
$$
If $\Vert h_0-h_1\Vert <\delta \sim\frac 1{\Vert f\Vert_{\rm Lip}}$, then
$$
\Vert \tau_{h_0} f - \tau_{h_1} f\Vert_2\leq (\Vert f\Vert_{\rm Lip} \delta)^{\frac 12} <\frac 12
$$
and consequently
$$
\Vert f-\tau_{h_1} f\Vert_2>1 \ \text { if } \ \Vert h_0-h_1\Vert<\delta.\eqno{(3.9)}
$$

In order to get a contradiction, we need to produce a word
$h_1 =g_1\cdots g_\ell; g_1, \ldots, g_\ell \in \mathcal G$ such that
$$
\Vert h_0 -g_1\cdots g_\ell\Vert <\delta\eqno{(3.10)}
$$
and
$$
\ell <\frac {\ve_1^4}{\sqrt\ve\sqrt d}.\eqno{(3.11)}
$$
Indeed, (3.6) implies then that
$$
\Vert f -\tau_{h_1} f\Vert_2 \leq \Vert f-\tau_{g_1} f \Vert_2 +\cdots + \Vert f-\tau_{g_\ell} f\Vert_2 < 1.
$$

For $1\leq i<d$, let $\sigma_{i, i+1} \in \text{Sym}(d)$ be the
transposition
 of $i$ and $i+1$.

Denote $\tilde \sigma_{i, i+1}$ the corresponding unitary operator. Since
$$
\{\sigma_{i, i+1} ; i =1, \ldots , d-1\}
$$
is a generating set for Sym$(d)$ consisting of cycles of bounded length, it follows from a result in
\cite{D-F} that the corresponding Cayley graph on
Sym$(d)$ has diameter at most $Cd^2$.
In particular, given $i, j \not\in \mathbb Z/d\mathbb Z, i\not= j, \tilde \sigma_{i, j}$ may be realized as a composition of a string of
elements $\tilde\sigma_{i, i+1}$ of length at most $Cd^2$.
In view of (3.7), this implies that if $\gamma\in SU(2)$ and $1\leq i<j\leq d$, then
$$
\Vert\gamma_{ij} -g\Vert < cd^2\ve_1\eqno{(3.12)}
$$
for some $g\in \mathcal G_{\ell_1}, \ell_1<cd^2$ ($\mathcal G_\ell$ = words of size $\ell$
written in $g$).

Let $\kappa>0$,
$$
\kappa^2> cd^2\ve_1.\eqno{(3.13)}
$$
Adopting the Lie-algebra point of view, the preceding implies that given
$s, t\in \mathbb R$, $|s|, |t|<1$ and $z\in\mathbb C$,
$|z|<1$, then
$$
\text{dist\,} \big(Id +\kappa\big( is(e_i\otimes e_i)+it (e_j\otimes e_j)+z(e_i\otimes e_j)+\bar z(e_j\otimes e_i)\big), \mathcal
G_{\ell_1}\big)<\kappa^2\eqno{(3.14)}
$$
and therefore
$$
\text{dist\,} (I+\kappa A, \mathcal G_{d^2\ell_1})<d^2 \kappa^2\eqno{(3.15)}
$$
for skew symmetric $A$, $\Vert A\Vert\leq 2\pi$.

Let $h\in SU(d), h= e^A$ with $A$ as above.
Taking $\kappa =\frac 1r$, we have
$$
e^A=(e^{\frac 1rA})^r =\Big(1+\frac 1r A\Big)^r +O\Big(\frac 1r\Big)
$$
and therefore, by (3.15)
$$
\text {dist\,} (h_0, \mathcal G_{rd^2\ell_1})\leq rd^2\kappa ^2 =\frac {d^2}r.\eqno{(3.16)}
$$
Taking $\kappa =\frac 1r =d^{-C}$ and $\ve_1=d^{-2C-2}$, (3.16) ensure that
$$
\text {dist\,} (h, \mathcal G_d c_1)< d^{-C} \ \text { for all } h\in SU(d).\eqno{(3.17)}
$$

Next, we rely on the Solovay-Kitaev commutator technique to produce approximations at smaller scale.
This procedure is in fact dimensional free (see the comment in \cite {D-N} following Lemma 2 in order to eliminate a polynomial prefactor
in $d$ - which actually would be harmless if we start from scales $\ve_0 =d^{-C}$).
The conclusion is that
$$
\rm{dist\,} (h, \mathcal G_\ell)< \tau \ \text { for all } \ h\in SU(d)
$$
may be achieved with
$$
\ell < d^{C_1} \Big(\log\frac 1\tau\Big)^A.
$$

Returning to (3.10), (3.11), we obtain the condition
$$
d^{C_0} \log^A(1+\Vert f\Vert_{\rm Lip}) <\frac {\ve_1^4}{\sqrt\ve\sqrt d} = d^{-C_2} \ve^{-\frac 12}\eqno{(3.19)}
$$
and Proposition 1$'$ follows.

\bigskip
\section
{Proof of Proposition 2}

The disadvantage of our approach is that $T$ is not restricted to finite dimensional invariant subspaces of $L^2(G)$
so that strictly speaking, one can not rely on a spectral gap argument to control the norm of iterates of $T$.

But Proposition 1$'$ nevertheless permit to derive easily the following

\begin{proposition}
Assume $f\in {\rm Lip}(G)$, $\Vert f\Vert_2 =1$, $\int_G f=0$.
Let $0<\rho<\frac 12$. Then
$$
\Vert T^\ell f\Vert_2<\rho\eqno{(4.1)}
$$
provided
$$
\ell> cd^C. \log^A (1+\Vert f\Vert_{\rm Lip}). \Big(\log \frac 1\rho\Big)^{A+1}.\eqno{(4.2)}
$$
\end{proposition}
\medskip

\begin{proof}
Let $B=\Vert f\Vert_{\rm Lip}$. Clearly $\Vert T^\ell f\Vert_{\rm Lip}\leq B$ also.

Fix some $\ell$ and let $f_1=\frac {T^\ell f}{\Vert T^\ell f\Vert_2}$.
Hence $\Vert f_1\Vert_{\rm Lip}\leq \frac B{\Vert T^\ell f\Vert_2}$.

Applying Proposition 1$'$, it follows that
$$
\Vert T^{\ell+1} f\Vert_2 \leq \Vert T^\ell f\Vert_2 (1-\ve_\ell)
$$
with
$$
\ve_\ell =cd^{-C}\Big(\log \Big(1+\frac B{\Vert T^\ell f\Vert_2}\Big)\Big)^{-A}>
cd^{-C} \big(\log(1+B)\big)^{-A} \Big(\log \Big(1+\frac 1{\Vert T^\ell f\Vert_2}\Big)\Big)^{-A}.
$$
Hence, assuming $\Vert T^\ell f\Vert_2>\rho$, we obtain
$$
\rho<(1-cd^{-C} \big(\log (1+B)\big)^{-A} \Big(\log\frac 1\rho\Big)^{-A}\Big)^\ell
$$
implying (4.2).
\end{proof}

\noindent
{\bf Proof of Proposition 2.}

Apply Proposition 3 with $\log B\sim \log \frac 1\ve$ and $\log \frac 1\rho\sim d^2 \log\frac 1\ve$.

\medskip

\section
{Variants}

The previous argument is clearly very flexible and may be applied in other situations.

Returning to \S3, assume more generally
$\eta$ a probability measure on $SU(2)$ satisfying
$\overline{\langle\supp \eta\rangle}= SU(2)$.
Note that by Proposition 1, $\eta^{(\ell)}$ with $\ell \sim(\log\frac 1{\ve_1})^c\sim (\log d)^c$ provides an
$\ve_1$-approximation of Haar measure on $SU(2)$.
It follows from (1.3), (3.3) that
$$
\int \Vert f-\tau_g f\Vert^2_2 (\nu_{i, i+1}*\nu_{i, i+1})(dg)< 2d^2\ve
$$
and hence
$$
\int\Vert f-\tau_g f\Vert^2_2 \nu^{(\ell)}_{i, i+1} (dg)<\ell d^2\ve
$$
$$
\nint _{\Omega_{a, i}} \Vert f-\tau_g f\Vert^2_2 \ \nu^{(\ell)}_{i, i+1}(dg)\lesssim \ell d^2 \ve_1^{-4}\ve.
$$
The collection $\mathcal G$ may then be introduced similarly.
Proposition 2 remains valid.

Let us point out that it is unknown if in general the density assumption
$\overline{\langle\supp\eta \rangle}=SU(2)$ implies a
uniform spectral gap (see the discussion on \S1).

Instead of (1.3), one may introduce at time $k=\mathbb Z_+$ the discrete average $T_k=\frac 12(\tau_g+\tau_{g-1})\cdots$
where we first pick
 some $i\in\mathbb Z/d\mathbb Z$ and then choose a random element $g\in
SU(2)$ acting on $[e_i,e_{i+1}]$ according to $\eta$.
In this situation, one obtains random walks on $SU(d)$ indexed by 
an additional probability space \break
$\otimes \big(\mathbb Z/d\mathbb Z \otimes SU(2)\big)$  
$$
T^\omega =\cdots T_k T_{k-1}\cdots T_1\eqno{(5.1)}
$$
and may ask for the typical mixing time of a realization.

Rather straightforward adjustments of the arguments appearing in the proof of Proposition $1'$ combined with some
Markovian considerations permit us to establish the analogue of Proposition 2 for $T^\omega$.
Thus

\begin{proposition}
Let $T^\omega$ be defined by (5.1).
Then, with large probability in $\omega$, $\ve$-approximation of uniform
measure on $SU(d)$ may be achieved in time $C(d
\log\frac 1\ve)^C$.
\end{proposition}

\noindent
{\bf Acknowledgement.}
The author is grateful to P.~Varju for his comments and also for pointing out various references.

P.~Varju also reported the following somewhat related question of A.~Lubotzky: Does $SU(d)$ admit a finite set of
generators with a spectral gap that is uniform in $d$?


\begin{thebibliography}
{XXXX}

\bibitem [B]{B} J.~Bourgain, \emph {A lower bound for the Lyapounov exponents of the random Schr\"odinger operator on a
strip}, J. Stat. Phys. 153 (2013), no 1, 1--9.

\bibitem [B-G1]{B-G1} J.~Bourgain, A.~Gamburd, \emph {On the spectral gap for finitely generated subgroups of $SU(2)$},
Inventiones Math A1 (2008), 201, 83--121.

\bibitem [B-G2]{B-G2} J.~Bourgain, A.~Gamburd, \emph{A spectral gap theorem in $SU(d)$}, JEMS 14 (2012), no5, 1455--1511.

\bibitem [B-L]{B-L} A.~Bougerol. J.~Lacroix, \emph {Product of random matrices with applications
to Schr\"odinger operators}, Birkh\"auser (1985).

\bibitem [B-S]{B-S} L.~Babai, A.~Seress: \emph {On the diameter of
permutation  groups}, European J.~Combin. 13 (1992),
231--243.

\bibitem[D-N]{D-N} C.~Dawson, M~Nielsen, \emph {The Solovay-Kitaev
algorithm},  Quantum Inf. \& Comp., Vol. 6, No 1 (2006) 81--95.

\bibitem [D-SC]{D-SC} P.~Diaconis, L.~Saloff-Coste, \emph{Bounds on the Kac's Master Equation}, CMP, 209 (3): 729--755
(2000).

\bibitem [D-F]{D-F} J.~Driscoll, M.~Furst, \emph{Computing short generator 
sequences}, Inf.  and Comp., Vol. 72, 1987, 117--132.

\bibitem[dS-B]{dS-B} Y.~Benoist, N.de Saxc\'e: \emph{ A spectral gap theorem  in simple Lie groups}, arXiv.1405.1808.

\bibitem[H-S]{H-S} H.~Helfgott, A.~Seress, \emph {On the diameter of permutation groups}, Annals Math. (2) 179 (2014), no
2, 611--658.

\bibitem [S-B]{S-B} H.~Schulz-Baldes, \emph {Perturbation theory for the Lyapounov exponents of 
an Anderson model on a strip}, GAFA 14, 1029-1117 (2004).


\bibitem [V]{V} P.~Varju, \emph{Random walks in compact groups}, Doc. Math. 18 (2013), 1137--1175.

\end{thebibliography}
\end{document}